\documentclass[11pt, a4paper, reqno, numberwithinsect,thm-restate,autoref]{amsart}
\usepackage[margin=1.2in]{geometry}
\usepackage{enumerate}

\usepackage{thmtools}
\usepackage{amsthm}

\usepackage{mathtools}

\numberwithin{equation}{section}
\numberwithin{table}{section}

\usepackage[colorlinks=true,linkcolor=blue]{hyperref}

\numberwithin{equation}{section}

\usepackage{amsmath}
\usepackage[mathscr]{euscript}
\usepackage{amssymb}
\usepackage[noadjust]{cite}
\usepackage{comment}
\usepackage{color}
\usepackage{filecontents}
\usepackage{bm}
\usepackage{aliascnt}

\makeatletter
\newcommand{\Eqref}[1]{\textup{\tagform@{\ref*{#1}}}}
\makeatother

\DeclareMathAlphabet{\mymathbb}{U}{bbold}{m}{n}

\usepackage[noabbrev]{cleveref}

\newcommand{\myref}[2][]{%
	\hyperref[#2]{\textup{(}\textnormal{\textit{#1}}\textup{)}}%
}

\newtheorem{thm}{Theorem}[section]
\newtheorem{lemma}[thm]{Lemma}

\newtheorem{prop}[thm]{Proposition}

\newtheorem{cor}[thm]{Corollary}

\theoremstyle{definition}

\theoremstyle{remark}

\numberwithin{equation}{section}

\newcommand{\BR}{\mathbb{R}}

\newcommand{\U}{\mathrm{U}}

\newcommand{\sphe}{\hbox{Sph}_{_{\RS_{\mathfrak{A}^+}}}}

\newcommand{\spheA}{\hbox{Sph}_{_{\RS_{{A}^+}}}}
\newcommand{\spheAJ}{\hbox{Sph}_{_{\RS_{{\mathfrak{A}}^+}}}}

\newcommand{\Proj}{\mathrm{Proj}}

\newcommand{\RS}{\mathrm{S}}

\newcommand{\one}{\mathbf{1}}

\newcommand{\half}{\frac{1}{2}}

\newcommand{\range}{r_{A^{**}}}

\numberwithin{equation}{section}

\begin{document}

\title[A metric characterization of projections among positive norm-one elements]{A metric characterization of projections among positive norm-One elements in unital C$^*$-algebras}

\author[A.M. Peralta and P. Saavedra]{Antonio M. Peralta \and Pedro Saavedra}

\address[A. M. Peralta]{Instituto de Matem{\'a}ticas de la Universidad de Granada (IMAG), Departamento de An{\'a}lisis Matem{\'a}tico, Facultad de
	Ciencias, Universidad de Granada, 18071 Granada, Spain.}
\email{aperalta@ugr.es}
\address[Pedro Saavedra]{Departamento de An\'{a}lisis Matem\'{a}tico, Facultad de Ciencias, Universidad de Granada, 18071 Granada, Spain.}
\email{psaavedraortiz@ugr.es}



\date{\today}

\begin{abstract} We characterize projections among positive norm-one elements in unital C$^*$-algebras in pure geometric terms determined by the norm of the underlying Banach space. Concretely, let $A$ be a C$^*$-algebra (or a JB$^*$-algebra) whose positive cone and unit sphere are denoted by ${A}^+$ and $\RS_{A}$, respectively. The positive portion of the unit sphere in $A$, denoted by $\RS_{{A}^+}$, is the set ${A}^+ \cap \RS_{A}$, while the unit sphere of positive norm-one elements around a subset $\mathscr{S}$ in $\RS_{A^+}$ is the set  $$\hbox{Sph}_{_{\RS_{{A}^+}}} (\mathscr{S}) :=\Big\{ x\in \RS_{{A}^+} : \|x-s\|=1 \hbox{ for all } s\in \mathscr{S} \Big\}.$$ Assuming that $A$ is unital, we establish that an element $a\in \RS_{{A}^+}$ is a projection if, and only if, it satisfies the double sphere property, that is,
	$ \hbox{Sph}_{_{\RS_{{A}^+}}} \left(\hbox{Sph}_{_{\RS_{{A}^+}}} \left(\{a\}\right) \right) = \{a\}.$ 
\end{abstract}
\maketitle

\section{Introduction}\label{sec:intro}

It is widely evidenced how the algebraic and geometric properties of C$^*$-algebras, and their elements, mutually interplay throughout deep and surprising interconnections. For instance, R.V. Kadison proved in 1951 that the set $\partial_e (\mathcal{B}_{A}),$ of extreme points of the closed unit ball $\mathcal{B}_{{A}}$ of a C$^*$-algebra $A$, is exactly the set of maximal partial isometries in $A$, that is, the elements $e\in A$ such that $ee^*$ (equivalently, $e^*e $) is a projection and $(\mathbf{1} - ee^*) A (\mathbf{1} - e^* e) = \{0\}.$ Moreover, in Kadison's terminology, the set of extreme points of the positive portion of the unit sphere in $A$ is the set of projections in $A$ (see \cite{Kad51}). A result by B. Aupetit showed that an element $a$ in a semi-simple Banach (or Jordan–Banach) algebra is an idempotent element if, and only if, its spectrum is contained in the set $\{0,1\}$ and there exist constants $r, C > 0$ such that $\sigma( x) \subset \sigma (a) + C \|x - a\|,$ for $\|x - a\| < r,$ where $\sigma (x)$ and $\sigma (a)$ stand for the spectrum of $x$ and $a$, respectively. C. Akemann and N. Weaver established in \cite{AkWe2002} a geometric characterization of partial isometries in C$^*$-algebras (see additionally \cite{FP2018}).\smallskip

Let $\mathscr{S}$ and $P$ be subsets of a Banach space $X$. The \emph{unit sphere around $\mathscr{S}$ in $P$} is defined as the set $$\hbox{Sph}_{_P}(\mathscr{S}) :=\Big\{ x\in P : \|x-b\|=1 \hbox{ for all } b\in \mathscr{S} \Big\}.$$ Note that these sets are determined by the norm on $X$. If we write ${A}^+$ and $\RS_{{A}^+}$ for the positive cone of $A$ and the positive portion of the unit sphere of a C$^*$-algebra $A$, respectively, we just need the norm to determine the set $\hbox{Sph}_{_{\RS_{{A}^+}}} (\mathscr{S}),$ for each $\mathscr{S}\subseteq A$. Given $a\in A$, we write $\hbox{Sph}_{_{\RS_{{A}^+}}}(a)$ for $\hbox{Sph}_{_{\RS_{{A}^+}}} (\{a\})$.\smallskip

Concerning the main goal of this note, it is established in \cite[Proposition 2.2]{peralta2018characterizing} that each positive norm-one element $a$ in a C$^*$-algebra $A$ satisfying $\hbox{Sph}_{_{\RS_{{A}^+}}} \left( \hbox{Sph}_{_{\RS_{{A}^+}}}(a)\right) =\{a\}$ is a projection. The question whether this metric condition is in fact a characterization of non-zero projections remains open. 
Positive answers to this challenging problem are known in the following particular cases: $A = M_n(\mathbb{C})$ by G. Nagy  \cite[Claim 1]{Nagy2018}, $A= B(H)$, where $H$ is a complex Hilbert space, by  the first author of this note \cite{peralta2018characterizing}, and $A$ being a von Neumann algebra or an AW$^*$-algebra by C.W. Leung, C.K. Ng, and N.C. Wong \cite{LNW25}. We extend the mentioned results by showing that the above question admits a positive solution when $A$ is a unital C$^*$-algebra (see \Cref{t chracterization of projections}).\smallskip

A \emph{JB$^*$-algebra} is a complex Banach space $\mathfrak{A}$ together with a bilinear and commutative product ($(a,b)\mapsto a\circ b$) and an algebra involution ($a\mapsto a^*$) satisfying the following axioms: 
\begin{enumerate}[$(1)$]
	\item $(a\circ b)\circ b^2  = (a\circ b^2)\circ b$ ($\forall a,b\in \mathfrak{A}$) (\emph{Jordan identity}),
	\item $\left\| \U_{a} (a^*) \right\| = \|a\|^3$ ($\forall a\in \mathfrak{A}$), where $\U_a (b) = 2 (a\circ b)\circ a - a^2 \circ b$  ($\forall a,b\in \mathfrak{A}$). 
\end{enumerate} A JB$^*$-algebra is called unital if there exists an element $\mathbf{1}\in \mathfrak{A}$ such that $\mathbf{1}\circ a = a\circ \mathbf{1} = a$ for all $a\in \mathfrak{A}$. Every C$^*$-algebra is a JB$^*$-algebra if we consider the Jordan product given by $a\circ b = \frac12 (a b + b a)$. However, the class of JB$^*$-algebras is strictly larger than the class of C$^*$-algebras. A JBW$^*$-algebra is a JB$^*$-algebra which is also a dual Banach space. As in the case of von Neumann algebras, JBW$^*$-algebras are always unital. Each JB$^*$-algebra $\mathfrak{A}$ admits a cone of positive elements, which is denoted by $\mathfrak{A}^+$. For the standard definitions and properties of JB$^*$-algebras employed in this paper, we refer to \cite{AlfsenShultz2003,HOS, Cabrera-Rodriguez-vol1} and to the detailed exposition in \cite[\S 2]{PeSa2025}.\smallskip 

In the Jordan setting the following properties hold for each  positive norm-one element $a$ in a JB$\,^*$-algebra $\mathfrak{A}$ (see \cite[Proposition 2.3]{PeSa2025}):  \begin{enumerate}[$(a)$]
	\item  	$\sphe(\sphe(a)) = \{a\}$ implies that $a$ is a projection in $\mathfrak{A}$.
	\item  If $\mathfrak{A}$ is a JBW$^*$-algebra, $a$ is a projection if, and only if, ${\sphe(\sphe(a)) = \{a\}}.$ 
\end{enumerate} The final section of this paper is devoted to proving that the previous equivalence in $(b)$ holds when $\mathfrak{A}$ is merely a unital JB$^*$-algebra. The final conclusion is that a positive norm-one element $p$ in a unital C$^*$-algebra (or in a unital JB$^*$-algebra) $A$ is a projection if, and only if, $\spheA\left(\spheA(p)\right) = \{p\}$.

\section{A metric characterization of projections in unital C$^*$-algebras}\label{Sec:Spheres Cstar}

The main result of this paper is the following metric characterization of projections among positive norm-one elements in a unital C$^*$-algebra.

%%% THEOREM 2.4
\begin{thm}\label{t chracterization of projections}
Let $A$ be a unital C$\,^*$-algebra, and let $p$ be a positive norm-one element in $A$. Then the following statements are equivalent:
	\begin{enumerate}[$(a)$]
		\item $p$ is a projection in $A$.
		\item $\spheA\left(\spheA(p)\right) = \{p\}$.
	\end{enumerate}
\end{thm} 

The proof of the previous theorem will be established after a series of technical results. Henceforth, the set of all projections in a C$^*$-algebra $A$ will be denoted by $\Proj(A)$. If $A$ is unital, the symbol $\mathbf{1}$ will stand for the unit element. Given $p\in \Proj(A)$ we shall write $p^{\perp}$ for $\mathbf{1}-p$. The spectrum of each $a\in A$ will be denoted by $\sigma(a).$  A well known property in C$^*$-algebra theory assures that for each norm-one element $x\in A$ and each projection $p\in A,$ the condition $pxp =p$ implies that $x = p + (\textbf{1}-p) x (\textbf{1}-p)$ (see, for example, \cite[Lemma 3.1]{FP2018}). We shall make us of this property without any further mention.\smallskip

The final paragraphs in the proof of $(a) \Rightarrow (b)$ in \cite[Theorem 2.3]{peralta2018characterizing} implicitly contain a weaker version of the following auxiliary result. We include here a stronger version with a simplified proof. 

%%% LEMMA 2.1
\begin{lemma}\label{lemma:1}
Let $A$ be a unital C$\,^*$-algebra. Then $\spheA\left(\spheA(\mathbf{1})\right) =\{\mathbf{1}\}$.
\end{lemma}

\begin{proof} Clearly $\mathbf{1}\in\spheA\left(\spheA(\mathbf{1})\right)$. Take $b\in \spheA\left(\spheA(\mathbf{1})\right)$. If $b$ is non-invertible (equivalently, $0$ lies in $\sigma(b)$), it can be easily seen that $\|b -\mathbf{1}\|=1,$ and hence $0= \|b-b\| =1$ by hypothesis, which is impossible. We can therefore assume that $b$ is invertible, or equivalently, $0\notin \sigma(b)$. If $\sigma(b) \neq \{1\}$, the number $m = \min \sigma(b)$ must be strictly smaller than $1$. It is routine to check that the element $\frac{b- m \mathbf{1}}{\|b- m \mathbf{1}\|}$ lies in $\RS_{{A}^+}$ with $\left\| \mathbf{1} - \frac{b- m \mathbf{1}}{\|b- m \mathbf{1}\|} \right\| = 1,$ and $\left\| b - \frac{b- m \mathbf{1}}{\|b- m \mathbf{1}\|} \right\| = m< 1,$ which contradicts that $b\in \spheA\left(\spheA(\mathbf{1})\right)$. Therefore $\sigma(b)=  \{1\},$ and thus $b= \mathbf{1}$.
\end{proof}

\begin{cor}\label{c original lemma 1} Let $p$ be a non-zero projection in a (probably non-unital) C$\,^*$-algebra $A$. Suppose  $b\in \spheA\left(\hbox{Sph}_{_{\RS_{(p {A} p)^+}}} (p)\right)$ with $b\leq p$. Then $b =p$.
\end{cor}

\begin{proof} Consider the unital C$^*$-algebra $p A p$, whose unit is precisely $p$. Since, by hypotheses, $b\in p A p$, it follows that $$b\in \spheA\left(\hbox{Sph}_{_{\RS_{(p {A} p)^+}}} (p)\right)\cap (p A p ) \subseteq  \hbox{Sph}_{_{\RS_{(p {A} p)^+}}} \left(\hbox{Sph}_{_{\RS_{(p {A} p)^+}}} (p)\right) = \{p\},$$ where the last equality follows from \Cref{lemma:1}. 
\end{proof}

Note that, under the assumptions in the previous lemma, the set  $\spheA\left(\hbox{Sph}_{_{\RS_{(p {A} p)^+}}} (p)\right)$ contains $\spheA\left(\spheA (p)\right)$.\smallskip

In the next lemma, we show that every element in the double positive sphere around a non-zero projection attains its norm at the Peirce-2 subspace associated with that projection. Since the last part of the note deals with JB$^*$-algebras, we shall try to use a terminology which is close to the Jordan setting, although no special Jordan theory will be required for those readers who are only interested in the result for C$^*$-algebras. Given an element $a$ in a C$^*$-algebra $A$, we write $\U_a$ for the operator on $A$ defined by $\U_a (b) = a b a$ ($b\in A$). Elements $a,b$ in $A$ are said to be orthogonal ($a\perp b$ in short) if $ a b^* = b^* a=0$. It is well-known that $a\perp b$ in $A$ implies $\|a+b\|= \max\{\|a\|,\|b\|\}.$\smallskip

Recall that a projection $p$ in a C$^*$-algebra $A$ is called minimal if it is non-zero and $p A p = \mathbb{C} p$.

%%% LEMMA 2.2
\begin{lemma}\label{lemma:2 pxp has norm one} Let $A$ be a unital C$\,^*$-algebra. Suppose $p$ is a non-zero projection in $A,$ and $b$ is an element in $\spheA\left(\spheA(p)\right)$. Then, $\|\U_p(b)\| = \| p b p \| = 1$.
\end{lemma}

\begin{proof} Set $c = \U_{\one-p}(b) = \U_{p^\perp}(b) \in \U_{p^\perp}(A) = p^\perp A p^\perp \subseteq A$. Observe that $c = 0$ if, and only if, $p^\perp \perp  b^{\half}$, equivalently, $p b = b p= b$, in which case $b = \U_p(b),$ and hence $\|\U_p(b)\| = \|b\| = 1$ as desired. We can therefore assume that $c \neq 0$.
	
Consider the continuous function $f : [0, \|c\|] \to [0, 1]$ given by $f(t) = \half + \frac{1}{2 \|c\|} t$ ($t \in [0, \|c\|]$).  Note that $f$ is affine and invertible, whose inverse is given by  $f^{-1}(s) = 2 \|c\| \left(s - \half\right)$ ($s \in [0, 1]$). Working on the C$^*$-subalgebra $ \U_{p^\perp}(A) ={p^\perp} A {p^\perp}$, we define, via continuous functional calculus, the element $d = f(c)\in {p^\perp} A {p^\perp}$ which clearly lies in $\mathrm{S}_{\U_{p^\perp}(A)}$. By construction $d \perp p$ and hence $\|d - p\| =\max \{\|d\|,\|p\|\} = 1,$ and thus $\|d - b\| = 1$ by hypotheses. Lemma 2.1 in \cite{peralta2018characterizing} assures the existence of a minimal projection $w \in A^{**}$ such that one of the following statements holds:
\begin{enumerate}[(I)]
		\item $w \leq d$ and $w \perp b$ in $A^{**}$.
		\item $w \leq b$ and $w \perp d$ in $A^{**}$.
	\end{enumerate}
	
\noindent Assume that statement (I) is true. Since $w,b \geq 0$, the condition $w \perp b$ is equivalent to 
	\begin{equation}\label{eq:wb_0}
		wb = bw = 0.
	\end{equation} It is well known that the inequalities $w \leq d \leq \mathbf{1}-p$ imply that $w\leq \mathbf{1}-p$ and $w = w b w$, and thus $$\begin{aligned}
	d = w + \U_{w^\perp}(d) &= w + (\mathbf{1}-w) b (\mathbf{1}-w) =  w + (\mathbf{1}-w) p^\perp b p^\perp (\mathbf{1}-w) \\
	&= w + (p^\perp-w) b (p^\perp-w) = w + \U_{(p^\perp-w)}(d) \in {p^\perp} A {p^\perp}.
\end{aligned}$$ A new application of the continuous functional calculus with the function $f^{-1}$ gives: 
\begin{equation}\label{eq:expression_of_c}
	c = f^{-1}\left(f(c)\right)= f^{-1}(d) = 2\|c\|\left(d - \half p^\perp\right) = 2\|c\|\left(w + \U_{(p^\perp-w)}(d) - \half p^\perp\right).
\end{equation} 

Having in mind that $w \leq p^\perp$, and hence $\U_w \circ \U_{p^\perp} = \U_w$, it follows from \eqref{eq:wb_0} and \eqref{eq:expression_of_c} that 
\begin{equation*}\begin{aligned}
		0 &= \U_w(b) = \U_w\left(\U_{p^\perp}(b) \right) = \U_w(c) =  2 \|c\| \U_w \left( w + \U_{(p^\perp-w)}(d) - \half p^\perp \right) \\
		&= 2 \|c\| \left( w^3 + \U_w(\U_{(p^\perp-w)}(d)) - \half \U_w(p^\perp) \right) = 2 \|c\| \left(w + 0 - \half w\right) = \|c\| w,
	\end{aligned}
\end{equation*} contradicting that $c$ and $w$ are both non-zero. Thus, case (I) above is impossible.\smallskip

\noindent Assume finally that statement (II) above holds. By construction $\frac12 p^\perp \leq d \leq p^\perp$, and hence the range projection of $d$ in $A^{**}$ (i.e. the smallest projection $q= \range(d)\in A^{**}$ satisfying $q d = d q = d$) coincides with $p^\perp$. Since $w \perp d$ by the assumptions in (II), we have $w \perp \range(d) = p^\perp$, and consequently $w \leq p$. It also follows from our assumptions that $w \leq b$. We therefore have $wp = pw = w = wb = bw$, and consequently $w = \U_w(b) = \U_w\left(\U_p(b)\right),$ which implies that
		\begin{equation*}
		1 = \|w\| = \|\U_w\left(\U_p(b)\right)\| \leq \|\U_p(b)\| \leq 1,
	\end{equation*} as desired.
\end{proof}

We can now get a refinement of our previous lemma. 

%%% PROPOSITION 2.3
\begin{prop}\label{prop:2.3 x b x has norm one} Let $A$ be a unital C$\,^*$-algebra. Suppose $p$ is a non-zero projection in $A,$ and $b$ is an element in $\spheA\left(\spheA(p)\right)$ with $\U_{p^\perp}(b)\neq 0$. Assume, additionally, that $x$ is an element in the set $\spheA(p) \cap (p A p) = \left\{ z\in p A^+ p : \|p - z\| = \|z\| = 1  \right\}$. Then, the element $\U_x(b) = x p x $ lies in the unit sphere of $A$.
\end{prop}

\begin{proof} Let us take $p$, $b$ and $x$ as in the statement. Set $y = p - x$, which is a positive element in the unit sphere of $p A p$ by assumptions. We consider first the following dichotomy:
	\begin{enumerate}[$(\bullet)$]
		\item If $\U_y(b) = 0$, we have $0 = yb = (p - x)b$, and so $pb = xb$. Taking adjoints we get $bp = bx$. We can therefore conclude that $\U_p(b) = \U_x(b)$ and hence, by Lemma \ref{lemma:2 pxp has norm one}, $\|\U_x(b)\| = \|\U_p(b)\| = 1$ as desired.
		\item If $\U_y(b) \neq 0$, we consider the element
		\begin{equation*}
			z = \half \left( p^\perp + \frac{ \U_{p^\perp}(b) }{ \| \U_{p^\perp}(b) \| } + \frac{ \U_{y}(b) }{ \| \U_{y}(b) \| } \right).
		\end{equation*}
	\end{enumerate} The next properties hold by construction: $0 \leq z$, $\left\| p^\perp + \frac{ \U_{p^\perp}(b) }{ \| \U_{p^\perp}(b) \| } \right\| = 2,$ $p^\perp + \frac{ \U_{p^\perp}(b) }{ \| \U_{p^\perp}(b) \|}\in  p^{\perp} A p^{\perp},$ and hence $\left(p^\perp + \frac{ \U_{p^\perp}(b) }{ \| \U_{p^\perp}(b) \| }\right) \perp \frac{ \U_{y}(b) }{ \| \U_{y}(b) \| }$. Therefore, by orthogonality, $$\|z\| = \max\left\{ \half \left\| p^\perp + \frac{ \U_{p^\perp}(b) }{ \| \U_{p^\perp}(b) \| } \right\|, \half \left\| \frac{ \U_{y}(b) }{ \| \U_{y}(b) \| } \right\|  \right\} = 1.$$ Similar considerations lead to $$\|z - p\| = \max\left\{ \half \left\| p^\perp + \frac{ \U_{p^\perp}(b) }{ \| \U_{p^\perp}(b) \| } \right\|, \half \left\| \frac{ \U_{y}(b) }{ \| \U_{y}(b) \| } - 2p \right\|  \right\} = 1.$$
	
The hypothesis on $b$ implies that $\|z - b\| = 1$. We can therefore apply \cite[Lemma 2.1]{peralta2018characterizing} to deduce the existence of a pure state $\varphi \in \partial_e\left(\mathcal{B}_{\left(A^*\right)^+ }\right)$ supported at a unique minimal projection $w = \mathrm{supp}(\varphi) \in A^{**}$ such that one of the following statements holds:
\begin{enumerate}[(I)]
\item $\varphi(z) = 1$ and $\varphi(b) = 0$.
\item $\varphi(z) = 0$ and $\varphi(b) = 1$.
\end{enumerate}
	
\noindent Assume that (I) holds. In such case we have $2 = 2 \, \varphi(z) = \varphi \left( p^\perp \right) + \varphi \left( \frac{ \U_{p^\perp}(b) }{ \| \U_{p^\perp}(b) \| } + \frac{ \U_{y}(b) }{ \| \U_{y}(b) \| } \right)$, with $\|p^\perp \|,$  $\left\| \frac{ \U_{p^\perp}(b) }{ \| \U_{p^\perp}(b) \| } + \frac{ \U_{y}(b) }{ \| \U_{y}(b) \| } \right\|\leq 1$. So, the only possibility is $$\varphi \left( p^\perp \right) = \varphi \left( \frac{ \U_{p^\perp}(b) }{ \| \U_{p^\perp}(b) \| } + \frac{ \U_{y}(b) }{ \| \U_{y}(b) \| } \right) = 1.$$ We consequently have 
	\begin{equation*}\label{eq:varphi_annihilates_p}
		\varphi(z) = \varphi\left({p^\perp}z\right)= \varphi\left(z {p^\perp} \right) = \varphi\left( {p^\perp} z {p^\perp} \right) = \varphi\left(\U_{p^\perp}(z)\right),  \text{ for all } z \in A.
	\end{equation*} In particular, $\varphi \left( \frac{ \U_{y}(b) }{ \| \U_{y}(b) \| } \right) = \frac{1 }{ \| \U_{y}(b) \| } \varphi \left( \U_{p^\perp} \left( \U_{y}(b)  \right)  \right) = \varphi (0) = 0$. Consequently, the equality $\varphi \left( \frac{ \U_{p^\perp}(b) }{ \| \U_{p^\perp}(b) \| } + \frac{ \U_{y}(b) }{ \| \U_{y}(b) \| } \right) = 1$ gives $\varphi \left(   \frac{ \U_{p^\perp}(b) }{ \| \U_{p^\perp}(b) \| } \right) = 1$. However, the condition $\varphi(b) = 0$ implies that
	\begin{equation*}
		0 = \varphi(b) = \varphi\left( \U_{p^\perp}(b)\right) = \| \U_{p^\perp}(b) \|,
	\end{equation*}
	which contradicts one of the hypotheses. Thus, case (I) is impossible.\smallskip
	
\noindent Suppose finally that case (II) above holds. The condition $\varphi(z) = 0,$ combined with the positivity of the elements $\U_y(b)$ and $p^\perp$, assures that $\varphi \left( \U_y(b) \right) = \varphi (p^\perp) = 0$, and the latter implies that $w \leq p$ (since $w$ is the minimal projection in $A^{**}$ supporting $\varphi$). Furthermore, by the Cauchy-Schwarz inequality, applied to the semi-positive sesquilinear mapping $(c, d) \mapsto \varphi(c d^*)$ at the elements $c = yb^\half$ and $d = b^\half$, we derive that
	\begin{equation}\label{eq:c-s}
		0 \leq | \varphi(yb) |^2 \leq \varphi(\U_y(b)) \; \varphi(b) = 0,
	\end{equation}
	so $\varphi(yb) = 0$ (and likewise $\varphi(by) = 0$).\smallskip 
	
	%Multiplying by the element $p - w \in A^{**}$, applying the Cauchy-Schwarz inequality once again and considering that $\varphi$ admits a {\color{red}canonical extension to $A^{**}$}, we get $\varphi \left( (p-w) b y\right) = 0$.
	Since $w$ is the support (minimal) projection of the pure state $\varphi$ in $A^{**}$, the identity $\U_w (a) = \varphi(a) w$ holds for all $a \in A$. In particular,
	\begin{equation*}\label{eq:wybyw_0}
		%\U_w \left( \U_y(b) \right) = \varphi \left( \U_y(b) \right) w = 0 \quad \text{and} \quad
		\U_w \left( yb \right) = \varphi \left( yb \right) w = 0.
	\end{equation*}
	To conclude the proof, we apply that $\varphi(b) = 1$ to rewrite the element $b$ in the form $b = w + \U_{\one - w}(b)$. In particular, $bw = wb = w$. This implies that $$0 = {\U_w(yb)  = w y b w = w y w = \U_w(y)} = \U_w(p - x),$$ which, in combination with $w \leq p$, yields $w = \U_w(p) = \U_w(x) = w x w$. It follows from the last identity that $x = w + w^{\perp} x w^{\perp}$, and thus $\U_w \circ \U_x = \U_w$. Finally, by considering the chain of inequalities given by
	\begin{equation*}
		1 \geq \|\U_x(b)\| \geq \| \U_w \left( \U_x(b) \right) \| = \| \U_w(b) \| = \|\varphi(b) \, w\| = 1,
	\end{equation*}
	we obtain the desired conclusion.
\end{proof}

\begin{proof}[Proof of Theorem~\ref{t chracterization of projections}]
	$(b) \Rightarrow (a)$. This is precisely Proposition 2.2 in \cite{peralta2018characterizing} (in fact, the result is true for arbitrary C$^*$-algebras).
	
	$(a) \Rightarrow (b)$ Let $p$ be a projection in $A$. It is clear that $\{p\} \subseteq \spheA\left(\spheA(p)\right)$. To prove the reverse inclusion, let us take an arbitrary $b \in \spheA\left(\spheA(p)\right)$. According to \Cref{c original lemma 1}, it suffices to prove one of the following equivalent conditions: $$b \leq p \Leftrightarrow \U_p(b) = p \Leftrightarrow \U_{p^\perp}(b) = 0.$$ Suppose, on the contrary, that $\U_{p^\perp}(b) \neq 0$.\smallskip
	
	Set $e = \U_p(b) \in \U_p(A)$. \Cref{lemma:2 pxp has norm one} implies that $\|e\| = 1$, and thus $e\in \RS_{p {A}^+ p}.$ Our goal will consist in showing that $e$ is the unit of the C$^*$-subalgebra $\U_p(A) = p A p$,  equivalently, $\sigma_{\U_p(A)} (e) = \{1\}$. The inclusion $\supseteq$ is clear because $\|e\| = 1$. To see the reverse inclusion, suppose, by contradiction, that there exists some value $t_0 \in \sigma_{\U_p(A)} (e) \cap [0, 1[$. 
	
	We can choose $0 < \varepsilon < \frac{1 - t_0}{2}$ and define, via continuous functional calculus in $p A p$, the element {$x = g(e)\in p A p$}, where $g : [0, 1] \to \BR$ is the continuous function given by
	\begin{equation*}
		g(t) = \begin{cases} 
			1, & \hbox{ if } 0 \leq t \leq t_0 + \varepsilon \\
			\frac{1 - \varepsilon - t}{1 - 2 \varepsilon - t_0}, & \hbox{ if } t_0 + \varepsilon \leq t \leq 1 - \varepsilon \\
			0, & \hbox{ if } 1 - \varepsilon \leq t \leq 1
		\end{cases}.
	\end{equation*}
	It is not difficult to see that $x \in \U_p(A) \cap \RS_{{A}^+}$ and $\|p - x\| = 1$. Since we assumed that $\U_{p^\perp}(b) \neq 0$, \Cref{prop:2.3 x b x has norm one} assures that $\|\U_x (b)\| = 1$. On the other hand, the fact that $x \in p A p$ implies that $\U_x(b) = \U_x(\U_p(b)) = \U_x(e)$, and thus
	\begin{equation*}
		1 = \|\U_x(b)\| = \|\U_x(e)\| = t_0 + \varepsilon < 1,
	\end{equation*}
	which is a contradiction. Therefore, no such $t_0$ exists, $\sigma_{\U_p(A)} (e) = \{1\},$ and $e = p$. This implies that $ p = \U_p (b) = p b p$, and consequently, $b = p + p^{\perp} b p^{\perp}$.\smallskip
	
Finally, we obtain a contradiction with the assumption that $z=\U_{p^\perp}(b)= p^{\perp} b p^{\perp}\neq 0.$ If $\|z\| =1$, we have $\| p - b\| = \|z\| =1$, equivalently, $b\in \spheA(p)$, and by hypothesis $\| b - b\| =1,$ which is impossible. Therefore $0<\|z\| <1$. For $0<\varepsilon <\|z\|$, consider the continuous function $h : [0, \|z\|] \to \BR$ given by
\begin{equation*}
	h(t) = \begin{cases} 
		t, & \hbox{ if } 0 \leq t \leq \|z\|- \varepsilon \\
		1 + \frac{1 - \|z\| + \varepsilon}{\varepsilon} (t - \|z\|), & \hbox{ if } \|z\|- \varepsilon \leq t \leq \|z\| \\
		1, & \hbox{ if } \|z\| \leq t \leq 1
	\end{cases},
\end{equation*} and the element $y = h(z)\in p^{\perp} A p^{\perp}$. It is easy to check that $0\leq y\leq p + y$ with  $\| p - (p+y)\| = \|y\|=1 = \|p +y\|$, and hence we must have $\| b - (p+y)\|=1$. However, $ \| b - (p-y)\|= \| p+ z - (p+y)\|=\|z-y\| = 1-\|z\| <1,$ which is impossible.  
\end{proof}

\section{The JB$^*$-algebra case}

Along this section we shall detail the necessary changes and alternative arguments to establish a generalization of our metric characterization of projections in \Cref{t chracterization of projections} to the wider setting of unital JB$^*$-algebras. We begin by observing that the norm-closed self-adjoint subalgebra of a unital JB$^*$-algebra generated by a self-adjoint element and the unit is isometrically Jordan $^*$-isomorphic to a commutative unital C$^*$-algebra (cf. \cite[Theorem 3.2.2 and \S 3.8]{HOS}). We can therefore conceive a continuous functional calculus for self-adjoint elements in a JB$^*$-algebra. Although we can not apply this property to deduce the next lemma from \Cref{lemma:1}, the argument can be easily adapted. 

\begin{lemma}\label{lemma:1 JBstar}
Let $\mathfrak{A}$ be a unital JB$\,^*$-algebra. Then $\spheAJ\left(\spheAJ(\mathbf{1})\right) =\{\mathbf{1}\}$.
\end{lemma}

\begin{proof} Let us take $b\in \spheAJ\left(\spheAJ(\mathbf{1})\right)$ and write $A$ for the JB$^*$-subalgebra of $\mathfrak{A}$ generated by $b$ and $\mathbf{1}$. As remarked above, $A$ is a commutative and unital C$^*$-algebra with unit $\mathbf{1}$. Furthermore, $A$ being generated by $b$ and $\mathbf{1}$ implies that $A\cong C(\sigma_{A} (b))$. The same argument in the proof of \Cref{lemma:1} shows that $\sigma_{A} (b) =\{1\}$, and hence $b=\mathbf{1},$ which concludes the proof.   
\end{proof}

The next corollary can be obtained through similar arguments to those applied in \Cref{c original lemma 1}. Note that, for each projection (i.e. a self-adjoint idempotent) $p$ in a JB$^*$-algebra $\mathfrak{A}$, the subspace $\U_p (A)$ is a JB$^*$-subalgebra of $\mathfrak{A}$ with unit $p$ (cf. \cite[\S 2.6]{HOS} or \cite[\S 1.3]{AlfsenShultz2003}).

\begin{cor}\label{c original lemma 1 JB} Let $p$ be a non-zero projection in a (probably non-unital) JB$\,^*$-algebra $\mathfrak{A}$. Suppose $b\in \spheAJ\left(\hbox{Sph}_{_{\RS_{\U_p (\mathfrak{A})^+}}} (p)\right)$ with $b\leq p$. Then $b =p$.
\end{cor}

The proof of \Cref{lemma:2 pxp has norm one} can be easily adapted to the Jordan setting up to suitable changes in the required tools. Let us recall that elements $a,b$ in a JB$^*$-algebra $\mathfrak{A}$ are orthogonal ($a\perp b$ in short) if $\{a,b,x\} = (a\circ b^*)\circ x +  (c\circ b^*)\circ a - (a\circ x)\circ b^* =0$ for all $x\in \mathfrak{A}$. In case that $a$ and $b$ are positive elements in $\mathfrak{A}$, we have $a\perp b \Leftrightarrow a\circ b=0$ \cite[Lemma 4.1]{BFGP09}, and in such a case $\U_a \U_b=0$. If we regard a C$^*$-algebra $A$ as a JB$^*$-algebra with respect to its natural Jordan product, the notion of orthogonality just defined agrees with the usual notion on C$^*$-algebras \cite[\S 2]{BFGMP08}. As in the associative setting, $a\perp b$ in $\mathfrak{A}$ implies $\|a+b\|= \max\{\|a\|,\|b\|\}.$ These properties together with several other reformulations of orthogonality are fully surveyed in \cite[\S 2]{BFGMP08} and \cite[\S 4]{BFGP09}.\smallskip

Another property required in our arguments writes as follows: \begin{equation}\label{eq Peirce is maximal}\left.\begin{aligned}
	&\hbox{Let $p$ be a projection in a JB$^*$-algebra $\mathfrak{A}$. Then, for each norm-one element}\\ &\hbox{ $x\in \mathfrak{A}$ with $\U_p(x) = p$ we have $x = p + \U_{p^{\perp}}(x)$ (see \cite[Lemma 1.6]{FriedmanRusso1985}).}
	\end{aligned}\right\}
\end{equation}

Given elements $a,b$ in a JB$^*$-algebra, the symbol $\U_{a,b}$ will stand for the linear mapping on $\mathfrak{A}$ defined by $\U_{a,b} (x) = (a\circ x)\circ b +  (b\circ x)\circ a - (a\circ b)\circ x.$ \smallskip

If $p$ is a projection in a JB$^*$-algebra $\mathfrak{A}$, the Jordan multiplication operator $M_{p} (x) = p\circ x$ ($x\in \mathfrak{A}$) has eigenvalues contained in the set $\{0,\frac12, 1\}$, the set $\mathfrak{A}_{j} (p) :=\{x\in \mathfrak{A}: p\circ x = \frac{j}{2} x \}$ is the eigenspace corresponding to the eigenvalue $\frac{j}{2}$ ($j\in \{0,1, 2\}$), and we have the following Peirce decomposition of $\mathfrak{A}$ with respect to $p$: $$ \mathfrak{A} = \mathfrak{A}_{0} (p)\oplus \mathfrak{A}_{1} (p)\oplus \mathfrak{A}_{2} (p).$$ The mappings $\U_p,$ $2\U_{p,\mathbf{1}-p}$, and $\U_{\mathbf{1}-p}$ are the projections of $\mathfrak{A}$ onto the direct summands $\mathfrak{A}_{2} (p)$, 
$\mathfrak{A}_{1} (p)$ and $\mathfrak{A}_{0} (p),$ respectively. The following multiplication rules hold: $$\begin{aligned}
\mathfrak{A}_{0} (p)\circ \mathfrak{A}_{0} (p) \subseteq \mathfrak{A}_{0} (p), \ \mathfrak{A}_{2} (p)\circ \mathfrak{A}_{2} (p) \subseteq \mathfrak{A}_{2} (p), \ \mathfrak{A}_{0} (p)\circ \mathfrak{A}_{2} (p) =\{0\},\ \ \ \ \ \ \ \ \ \ \\ \left(\mathfrak{A}_{0} (p)\oplus  \mathfrak{A}_{2} (p)\right) \circ \mathfrak{A}_{1} (p) \subseteq \mathfrak{A}_{1} (p), \hbox{ and } \mathfrak{A}_{1} (p)\circ \mathfrak{A}_{1} (p) \subseteq \left(\mathfrak{A}_{0} (p)\oplus  \mathfrak{A}_{2} (p)\right) \ \hbox{(see \cite[\S 2.6]{HOS})}. 
\end{aligned} $$ It follows from this rules that 
$$\U_{\mathfrak{A}_{0} (p),\mathfrak{A}_{2} (p)} \left(\mathfrak{A} \right) = \U_{\mathfrak{A}_{0} (p),\mathfrak{A}_{2} (p)} \left(\mathfrak{A}_{1} (p) \right) \subseteq \mathfrak{A}_{1} (p),$$ 
$$\U_{\mathfrak{A}_{0} (p),\mathfrak{A}_{0} (p)} \left(\mathfrak{A}_{1} (p) \right)\subseteq \mathfrak{A}_{1} (p), \hbox{ and hence } p\circ \U_{\mathfrak{A}_{0} (p),\mathfrak{A}_{0} (p)} \left(\mathfrak{A}_{1} (p) \right) = \frac12 \U_{\mathfrak{A}_{0} (p),\mathfrak{A}_{0} (p)} \left(\mathfrak{A}_{1} (p) \right).$$ However, given $a_0,b_0\in \mathfrak{A}_{0} (p)$ and $c_1\in \mathfrak{A}_{1} (p),$ by \cite[Lemma 2.4.20]{HOS} we have
$$ \begin{aligned}
\frac12 \U_{a_{0},b_{0}} \left(c_{1} \right) &= p\circ \U_{a_{0},b_{0}} \left(c_{1} \right) = \U_{p\circ a_{0} , b_{0}} \left(c_{1} \right) + \U_{a_{0}, p\circ  b_{0}} \left(c_{1} \right)- \U_{a_{0}, b_{0}} \left(p\circ   c_{1} \right) = - \frac12 \U_{a_{0},b_{0}} \left(c_{1} \right),
\end{aligned} $$ which guarantees that \begin{equation}\label{eq Peirce 010} \U_{\mathfrak{A}_{0} (p),\mathfrak{A}_{0} (p)} \left(\mathfrak{A}_{1} (p) \right) =\{0\}, \hbox{ and consequently } \U_{\mathfrak{A}_{0} (p),\mathfrak{A}_{0} (p)} \left( \U_{\mathfrak{A}_{0} (p),\mathfrak{A}_{2} (p)} \left(\mathfrak{A} \right)   \right) =\{0\}.
\end{equation}

As in the setting of C$^*$-algebras, a projection $p$ in a JB$^*$-algebra $\mathfrak{A}$ is called minimal if $\U_{p} (\mathfrak{A}) = \mathbb{C} p\neq \{0\}.$ A state on a JB$^*$-algebra (i.e. a positive norm-one functional) is pure if it is not a convex combination of two distinct states. For each pure state $\varphi$ on a JB$^*$-algebra $\mathfrak{A}$ there exists a unique minimal projection $w \in \mathfrak{A}^{**}$ satisfying $\U_w (x) = \varphi(x) w$ for all $x\in \mathfrak{A}$ (see \cite[Lemma 4.14]{AlfsenShultz2003} or \cite[Proposition 4]{FriedmanRusso1985} for a more general result).

\begin{lemma}\label{lemma:2 pxp has norm one JB} Let $\mathfrak{A}$ be a unital JB$\,^*$-algebra. Suppose $p$ is a non-zero projection in $\mathfrak{A},$ and $b$ is an element in $\spheAJ\left(\spheAJ(p)\right)$. Then, $\|\U_p (b)\| = 1$.
\end{lemma}

\begin{proof} Define $c := \U_{\one-p}(b) = \U_{p^\perp}(b) \in \U_{p^\perp}(\mathfrak{A}) \subseteq \mathfrak{A}$. Observe that $c = 0$ if, and only if, $p^\perp \circ b^{\half} = 0,$ equivalently, $p^\perp \circ b =0\Leftrightarrow p\circ  b = b$, in which case $b = \U_p(b),$ and hence $\|\U_p(b)\| = \|b\| = 1$ as desired. There is no loss of generality in assuming that $c \neq 0$. As in the associative setting, we define, via continuous functional calculus in the unital JB$^*$-algebra $\U_{p^\perp}(\mathfrak{A})$, the element $d = f(c)\in \U_{p^\perp}(\mathfrak{A}),$ where $f : [0, \|c\|] \to [0, 1]$ is given by $f(t) = \half + \frac{1}{2 \|c\|} t$ ($t$ in $[0, \|c\|]$). It is easy to check that, by orthogonality, $\|d - p\| =\max \{\|d\|,\|p\|\} = 1,$ and thus $\|d - b\| = 1$ by hypotheses. By \cite[Lemma 2.1$(a)$]{PeSa2025}, \eqref{eq Peirce is maximal}, and the comments before this lemma, there exists a minimal projection $w \in \mathfrak{A}^{**}$ such that one of the following statements holds:
	\begin{enumerate}[(I)]
		\item $d = w + \U_{w^\perp}(d)$ and $w \perp b$ in $\mathfrak{A}^{**}$.
		\item $b = w + \U_{w^\perp}(b)$ and $w \perp d$ in $\mathfrak{A}^{**}$.
	\end{enumerate}
We shall consider both cases independently.\smallskip

 \noindent(I) Since $w \leq d \leq p^\perp,$ we have $p^\perp = w + (p^\perp -w),$ with $w \perp (p^\perp - w)$. We can therefore deduce that $$\begin{aligned}
	d = w + \U_{w^\perp}(d) &= w + \U_{w^\perp} \U_{p^\perp} (d) = w + \U_{w^\perp} \left(\U_{w} (d) + 2 \U_{w, p^\perp-w} (d) + \U_{p^\perp-w} (d) \right) \\
	&= w + \U_{w^\perp} \U_{w} (d) + 2  \U_{w^\perp} \U_{w, p^\perp-w} (d) +  \U_{w^\perp} \U_{p^\perp-w} (d)   \\
	&= w +  2  \U_{w^\perp} \U_{w, p^\perp-w} (d) +  \U_{w^\perp} \U_{p^\perp-w} (d) .
\end{aligned}$$ It follows from \eqref{eq Peirce 010} that $\U_{w^\perp} \U_{w, p^\perp-w} (d) =0.$ Moreover, since  $w^\perp \circ (p^\perp-w) = p^\perp-w,$ and hence $\U_{p^\perp-w} (d) \in \mathfrak{A}_{2} (w^\perp)$, we can assert that $\U_{w^\perp} \U_{p^\perp-w} (d) = \U_{p^\perp-w} (d) $. Thus, $d = w + \U_{p^\perp-w} (d).$ A new application of the continuous functional calculus in  $\U_{p^\perp}(\mathfrak{A})$ now gives: 
\begin{equation}\label{eq:expression_of_c Jordan}
c = f^{-1}\left(f(c)\right)= f^{-1}(d) = 2\|c\|\left(d - \half p^\perp\right) = 2\|c\|\left(w + \U_{(p^\perp-w)}(d) - \half p^\perp\right).
\end{equation} On the other hand, having in mind that $w\perp b,$ $w\leq p^\perp$ and the previous conclusions, we get 
\begin{equation*}\begin{aligned}
		0 &= \U_w(b) = \U_w\left(\U_{p^\perp}(b) \right) = \U_w(c) =  2 \|c\| \U_w \left( w + \U_{(p^\perp-w)}(d) - \half p^\perp \right) \\
		&= 2 \|c\| \left( w^3 + \U_w(\U_{(p^\perp-w)}(d)) - \half \U_w(p^\perp) \right) = 2 \|c\| \left(w + 0 - \half w\right) = \|c\| w \neq 0,
	\end{aligned}
\end{equation*}  which is a contradiction. So, case (I) above is impossible.\smallskip

\noindent (II) By construction $\frac12 p^\perp \leq d \leq p^\perp$, and hence the range projection of $d$ in $\mathfrak{A}^{**}$ (cf. \cite[Lemma 4.2.6]{HOS}), coincides with $p^\perp$. Since $w \perp d$ by the assumptions in (II), we have $w \perp \range(d) = p^\perp$, and consequently $w \leq p$. It also follows from our assumptions that $b = w + \U_{w^\perp}(b),$ and thus $\U_w(b) = \U_w(w + \U_{w^\perp}(b)) = w,$ which implies that
	\begin{equation*}
		1 = \|w\|  = \|\U_w\left(b\right)\| = \|\U_w\left(\U_p(b)\right)\| \leq \|\U_p(b)\| \leq 1,
	\end{equation*} as desired.
\end{proof}

In order to establish a Jordan version of \Cref{prop:2.3 x b x has norm one} we shall need a refinement of the Cauchy-Schwarz inequality.  

\begin{lemma}\label{l CS inequality refined} Let $\varphi$ be a positive functional on a JB$\,^*$-algebra $\mathfrak{A}$. Suppose $a$ and $b$ are positive elements in $\mathfrak{A}$ such that $\varphi(\U_a (b)) = 0$. Then $\varphi(a \circ b) = 0$.
\end{lemma}

\begin{proof} There is no problem in assuming that $\mathfrak{A}$ is unital with unit $\mathbf{1}$. Let $B$ denote the JB$^*$-subalgebra of $\mathfrak{A}$ generated by $a,b$ and $\mathbf{1}$. By the Shirshov-Cohn theorem (cf. \cite[Theorems 2.4.14 and 7.2.5]{HOS} or \cite[Corollary 2.2]{Wright1977}) there exists a unital C$^*$-algebra $A$ with unit $\mathbf{1}$ containing $B$ as a JB$^*$-subalgebra. The product on $A$ will be denoted by mere juxtaposition. Clearly, $\varphi|_{B}$ is a positive functional in $B^*$. Let $\phi$ denote a Hahn-Banach extension of $\varphi|_{B}$ to $A$. We can see that $\phi$ is a positive functional on $A$ \cite[Proposition 1.5.2]{Sak}. Furthermore, $0=\varphi(\U_a (b)) = \phi(a b a)$. Finally, by the Cauchy-Schwarz inequality we derive that 
	$$ |\varphi(a \circ b)| = |\phi(a \circ b)| \leq  \frac12 \left( |\phi(a b^{\frac12} b^{\frac12} )| + |\phi( bb^{\frac12} b^{\frac12} a )|  \right)\leq \phi (b)^{\frac12} \phi( aba )^{\frac12} = 0,$$ which concludes the proof.
\end{proof}

The Jordan version of \Cref{prop:2.3 x b x has norm one} now reads as follows: 
 
\begin{prop}\label{prop:2.3 x b x has norm one Jordan} Let $\mathfrak{A}$ be a unital JB$\,^*$-algebra. Suppose $p$ is a non-zero projection in $\mathfrak{A},$ and $b$ is an element in $\spheAJ\left(\spheAJ(p)\right)$ with $\U_{p^\perp}(b)\neq 0$. Let us assume that $x$ is an element in the set $\spheAJ(p) \cap \U_p \left(\mathfrak{A}\right) = \left\{ z\in \U_p \left(\mathfrak{A}^+\right) : \|p - z\| = \|z\| = 1  \right\}$. Then, the element $\U_x(b)$ lies in the unit sphere of $\mathfrak{A}$. 
\end{prop}

\begin{proof} The proof is very similar to that of \Cref{prop:2.3 x b x has norm one}, and we shall try to simplify the details. Take $y = p - x\in \RS_{{\mathfrak{A}}^+}\cap \U_{p} (\mathfrak{A})$. If $\U_y(b) = 0$ (equivalently, $y\perp b$), by the Shirshov-Cohn theorem (cf. \cite[Theorems 2.4.14 and 7.2.5]{HOS} or \cite[Corollary 2.2]{Wright1977}), denoting by $B$ the JB$^*$-subalgebra of $\mathfrak{A}$ generated by $y,b$ and $\mathbf{1}$, there exists a unital C$^*$-algebra $A$ with unit $\mathbf{1}$ containing $B$ as a JB$^*$-subalgebra. Employing the (associative) product of $A$, which will be denoted by juxtaposition, we have $p b = x b,$ $b p = b x$, and $\U_p (b) = p b p = x b x= \U_x(b)$. In such a case Lemma \ref{lemma:2 pxp has norm one JB} assures that $\|\U_x(b)\| = \|\U_p(b)\| = 1,$ as desired. We can therefore assume that  $\U_y(b) \neq 0$.\smallskip 
	
Consider the element
		\begin{equation*}
			z = \half \left( p^\perp + \frac{ \U_{p^\perp}(b) }{ \| \U_{p^\perp}(b) \| } + \frac{ \U_{y}(b) }{ \| \U_{y}(b) \| } \right).
		\end{equation*} It is easy to check that $0 \leq z$, $\left\| p^\perp + \frac{ \U_{p^\perp}(b) }{ \| \U_{p^\perp}(b) \| } \right\| = 2,$ $p^\perp + \frac{ \U_{p^\perp}(b) }{ \| \U_{p^\perp}(b) \|}\in  \U_{p^{\perp}} \left(\mathfrak{A}\right)^{+},$ with $\left(p^\perp + \frac{ \U_{p^\perp}(b) }{ \| \U_{p^\perp}(b) \| }\right) \perp \frac{ \U_{y}(b) }{ \| \U_{y}(b) \| }$, $$\|z\| = \max\left\{ \half \left\| p^\perp + \frac{ \U_{p^\perp}(b) }{ \| \U_{p^\perp}(b) \| } \right\|, \half \left\| \frac{ \U_{y}(b) }{ \| \U_{y}(b) \| } \right\|  \right\} = 1,$$ and  $$\|z - p\| = \max\left\{ \half \left\| p^\perp + \frac{ \U_{p^\perp}(b) }{ \| \U_{p^\perp}(b) \| } \right\|, \half \left\| \frac{ \U_{y}(b) }{ \| \U_{y}(b) \| } - 2p \right\|  \right\} = 1.$$ Therefore, the hypothesis on $b$ implies that $\|z - b\| = 1$. Lemma 2.1 in \cite{PeSa2025} implies the existence of a pure state $\varphi$  of $\mathfrak{A}$ supported at a unique minimal projection $w = \mathrm{supp}(\varphi)$ in $\mathfrak{A}^{**}$ such that one of the following statements holds:
	\begin{enumerate}[(I)]
		\item $\varphi(z) = 1$ and $\varphi(b) = 0$.
		\item $\varphi(z) = 0$ and $\varphi(b) = 1$.
	\end{enumerate}
	
\noindent In case (I) we have $2 = 2 \, \varphi(z) = \varphi \left( p^\perp \right) + \varphi \left( \frac{ \U_{p^\perp}(b) }{ \| \U_{p^\perp}(b) \| } + \frac{ \U_{y}(b) }{ \| \U_{y}(b) \| } \right)$, and hence $$\varphi \left( p^\perp \right) = \varphi \left( \frac{ \U_{p^\perp}(b) }{ \| \U_{p^\perp}(b) \| } + \frac{ \U_{y}(b) }{ \| \U_{y}(b) \| } \right) = 1.$$ In particular $p^{\perp}\geq w,$ $\varphi(x) = \varphi\left(\U_{p^\perp}(z)\right)\,$  for all $z \in \mathfrak{A}$, and thus $\varphi \left( \frac{ \U_{y}(b) }{ \| \U_{y}(b) \| } \right) = \frac{1 }{ \| \U_{y}(b) \| } \varphi \left( \U_{p^\perp} \left( \U_{y}(b)  \right)  \right) = \varphi (0) = 0$. Consequently, $\varphi \left(   \frac{ \U_{p^\perp}(b) }{ \| \U_{p^\perp}(b) \| } \right) = 1$. However, the fact that $\varphi(b) = 0$ implies that
	\begin{equation*}
		0 = \varphi(b) = \varphi\left( \U_{p^\perp}(b)\right) = \| \U_{p^\perp}(b) \|\neq 0,
	\end{equation*}
	which is impossible.\smallskip
	
\noindent (II) The condition $\varphi(z) = 0,$ combined with the positivity of the elements $\U_y(b)$ and $p^\perp$, assures that $\varphi \left( \U_y(b) \right) = \varphi (p^\perp) = 0$, and the latter implies that $w \leq p$ (since $w$ is the minimal projection in $\mathfrak{A}^{**}$ supporting $\varphi$). \Cref{l CS inequality refined} assures that $\varphi (y\circ b) =0$. Since $w$ is the support (minimal) projection of the pure state $\varphi$ in $A^{**}$, the identity $\U_w (a) = \varphi(a) w$ holds for all $a \in A$. In particular,
\begin{equation}\label{eq: Jan 14 1}
	\U_w \left( y \circ b \right) = \varphi \left( y\circ b \right) w = 0.
\end{equation} The condition $\varphi(b) = 1$ and the comments prior to \Cref{lemma:2 pxp has norm one JB} give $b = w + \U_{w^{\perp}}(b)$, and by \eqref{eq: Jan 14 1}, \begin{equation}\label{eq Jan 14 b} 0 = \U_{w} (y\circ b) = \U_{w} (y\circ w) + \U_{w} (y\circ \U_{w^{\perp}}(b)).
\end{equation} By applying the Peirce decomposition of $\mathfrak{A}^{**}$ with respect to $w$, we deduce that the element $y\circ \U_{w^{\perp}}(b)$ lies in the space $  \mathfrak{A}^{**}_{0} (w) \oplus \mathfrak{A}^{**}_1 (w),$ and hence $ \U_{w} (y\circ \U_{w^{\perp}}(b)) =0$. Another application of the Peirce decomposition shows that $\U_{w} (y\circ w) = \U_{w} (y).$ It follows from \eqref{eq Jan 14 b} that $\U_{w} (y)=0$, equivalently, $\U_w (x) = \U_w (p) =w$. An application of \eqref{eq Peirce is maximal} proves that $ x = w +  \U_{w^{\perp}}(x)$. It follows from the above conclusions and the properties of the Peirce decomposition that 
$$ \U_{w, {\U_{w^{\perp}}(x)}} (b) = \U_{w, \U_{w^{\perp}}(x)} (w + \U_{w^{\perp}}(b)) = 0 ,$$ and 
$\U_{\U_{w^{\perp}}(x)} (b) \in \mathfrak{A}^{**}_0 (w).$ Finally, the inequalities
$$\begin{aligned}
1 \geq \|\U_x(b)\| &\geq \| \U_w (b) + 2 \U_{w, \U_{w^{\perp}}(x)} (b) + \U_{\U_{w^{\perp}}(x)} (b) \| \\
&=  \| \U_w (b) + \U_{\U_{w^{\perp}}(x)} (b) \| = \max\{ \| \U_w(b) \|, \| \U_{\U_{w^{\perp}}(x)} (b) \|\} \\
&\geq \| \U_w(b) \| = \|\varphi(b) \, w\| = 1
\end{aligned}$$ give the desired result.
\end{proof}

We can now obtain a version of \Cref{t chracterization of projections} for unital JB$^*$-algebras. 

\begin{thm}\label{t chracterization of projections JB}
	Let $\mathfrak{A}$ be a unital JB$\,^*$-algebra, and let $p$ be a positive norm-one element in $\mathfrak{A}$. Then the following statements are equivalent:
	\begin{enumerate}[$(a)$]
		\item $p$ is a projection in $\mathfrak{A}$.
		\item $\spheAJ\left(\spheAJ(p)\right) = \{p\}$.
	\end{enumerate}
\end{thm}

\begin{proof} The implication $(b)\Rightarrow (a)$ is proved in \cite[Proposition 2.3$(a)$]{PeSa2025}. Conversely, when in the proof of \Cref{t chracterization of projections}$(a)\Rightarrow (b)$ we replace \Cref{c original lemma 1}, \Cref{lemma:2 pxp has norm one} and \Cref{prop:2.3 x b x has norm one} with \Cref{c original lemma 1 JB}, \Cref{lemma:2 pxp has norm one JB}, and \Cref{prop:2.3 x b x has norm one Jordan}, respectively, the argument remains valid. 
\end{proof}

\medskip

\section*{Acknowledgement}

A. M. Peralta supported by grant PID2021-122126NB-C31 funded by MICIU/AEI/ 10.13039/501100011033 and by ERDF/EU, by Junta de Andalucía grant FQM375, IMAG--Mar{\'i}a de Maeztu grant CEX2020-001105-M/AEI/10.13039/501100011033 and (MOST) Ministry of Science and Technology of China grant G2023125007L. \smallskip

P. Saavedra supported by a Formaci\'{o}n de Profesorado Universitario (FPU) grant from the Ministerio de Ciencia, Innovaci\'{o}n y Universidades, Spain (Grant No. FPU23/00747). \medskip

\smallskip\smallskip

\noindent\textbf{Data Availability} Data sharing is not applicable to this article as no datasets were generated or analysed during the preparation of the paper.\smallskip\smallskip

\noindent\textbf{Declarations} 
\smallskip\smallskip

\noindent\textbf{Conflict of interest} The authors declare that they have no conflict of interest.

\end{document}